\documentclass{amsart}

\usepackage{amsmath,amsfonts,amssymb,amsthm,amscd,latexsym,euscript,xypic}
\usepackage[all]{xy}

\newtheorem {theorem1}{Theorem}[section]
\newtheorem {theorem}[theorem1]{Theorem}

\newtheorem {lemma}[theorem1]{Lemma}
\theoremstyle{definition}

\theoremstyle{remark}

\newcommand{\calA}{\ensuremath{\mathcal{A}}}

\newcommand{\calD}{\ensuremath{\mathcal{D}}}
\newcommand{\calS}{\ensuremath{\mathcal{S}}}
\newcommand{\calX}{\ensuremath{\mathcal{X}}}

\newcommand{\cat}[1]{\ensuremath{\EuScript #1}}

\DeclareMathOperator{\map}{\textup{map}}

\newcommand{\colim}{\ensuremath{\mathop{\textup{colim}}}}

\newcommand{\rarrow}{\rightarrow}
\newcommand{\Id}{\ensuremath{\textup{Id}}}
\newcommand{\op}{\ensuremath{\textup{op}}}

\newcommand{\Vopenkas}{Vo\-p\v en\-ka's }
\newcommand{\Rosicky}{Rosick{\'y} }

\newcommand{\hbot}{\ensuremath{\mathop{{}^\mathrm{h}\!\!\bot}}}
\newcommand{\htop}{\ensuremath{\mathop{{}_\mathrm{h}\!\!\top}}}

\newcommand{\CW}{\ensuremath{{\textup{CW}}}}

\begin{document}

\SelectTips{cm}{10}

\title[Cellularization with respect to a set of objects]
      {Abstract cellularization as a cellularization with respect to a set of objects}
\author{Boris Chorny}
\thanks{}

\address{Department of Mathematics, ETH Zentrum, 8092 Zurich, Switzerland}

\email{chorny@math.ethz.ch}

\date{\today}
\dedicatory{} \commby{}

\begin{abstract}
Given a simplicial idempotent augmented endofunctor $F$ on a simplicial combinatorial model category \cat M, under the assumption of \Vopenkas principle, we exhibit a set $A$ of cofibrant objects in \cat M such that $F$ is equivalent to $\CW_A$, the cellularization with respect to $A$.
\end{abstract}

\maketitle
\section*{Introduction}
Over the past ten years many results about locally presentable
categories were generalized to combinatorial model categories.
The term ``combinatorial'' means that the model category is cofibrantly
generated and the underlying category is locally presentable (see
\cite{AR} and \cite{Hirschhorn} for the definitions of these
concepts). This notion is due to J.~H.~Smith,  who constructed (in
unpublished work) localizations of combinatorial model category
structures with respect to sets of maps. This construction may be
viewed as a generalization of the orthogonal reflection in a locally
presentable model category \cite[1.35]{AR} and the claim that every
small-orthogonality class of a locally presentable category is locally
presentable \cite[1.40]{AR}. In \cite{Dugger:presentations}, Dugger
proved that every combinatorial model category is equivalent to a
localization of a category of diagrams of simplicial sets, hence
generalizing \cite[1.46]{AR}. Furthermore, J.~\Rosicky have recently
proved (personal communication) that for strongly combinatorial model
category \cat K, (cone) injectivity classes in the homotopy category
of \cat K behave similarly as in locally presentable
categories. \cat K is \emph{strongly combinatorial} if it is
combinatorial and the class of cofibrations is closed under $\kappa$-directed colimits in $\mathop{\mathrm{Map}}(\cat K)$
for some regular cardinal $\kappa$. 

\Vopenkas principle has especially powerful implications for locally
presentable categories if one is ready to assume it
\cite[Chapter~6]{AR}. In particular existence of localization with
respect to an arbitrary class of maps in a locally presentable
category is equivalent to the \textbf{weak} \Vopenkas principle. Substituting
weak \Vopenkas principle by a stronger \Vopenkas principle, J.\Rosicky
and W.~Tholen \cite{Rosicky-Tholen} proved an analogous result for
combinatorial model categories: \Vopenkas principle implies the
existence of homotopical localizations in combinatorial model
categories with respect to arbitrary classes of maps. In the category
of simplicial sets this result was proven earlier in the work of
C.~Casacuberta, D.~Scevenels, and J.~Smith \cite{CSS}.

If \Vopenkas principle is assumed, a
stronger statement for locally presentable categories holds: any
orthogonality class is a small orthogonality class (i.e., it is generated by a \textbf{set} of maps). An analog for
combinatorial model categories is the following: any homotopy
localization functor is a localization with respect to a set of maps
\cite{cacho}. There exists a counterexample for this statement in a
locally presentable category which is not cofibrantly generated
\cite{PhDII}, so combinatorial model categories are the right
generalization of locally presentable model categories.

We will finish this comparison of locally presentable categories with
combinatorial model categories with one example of a statement which was
not generalized yet: any orthogonality class generated by a set of
objects in a complete and well-powered category was shown to be
reflective in \cite{CFT}. An analogous statement in topology is known better
as a problem of existance of homotopical localization with respect to
a cohomology theory. Cohomological localizations are known to exist
under \Vopenkas principle, but in ZFC this problem is still open.

In this paper we continue the efforts of transferring results from
locally presentable category theory to abstract homotopy
theory. Namely we generalize the following statement: assuming
\Vopenkas principle any co-orthogonality class in a locally
presentable category is co-reflective \cite[6.29]{AR}. We show that
under \Vopenkas principle any homotopy co-orthogonality class in a
combinatorial model category  may be
generated by a set of objects and hence (if the category is right proper) admits a homotopy
co-localization (e.g., by \cite[Prop. 5.3.5]{Hirschhorn} applied to a map
$\emptyset\rarrow X$).

Next we point out that any cellularization (homotopy co-localization)
in a simplicial combinatorial model category gives rise to a homotopy co-orthogonal pair, and hence
it is equivalent to  a cellularization with respect to a set
of objects. This question is not formally dual to the one discussed in
\cite{cacho}, since the opposite of a locally presentable category is
usually no longer locally presentable.

If one considers an arbitrary locally presentable category as a model
category with weak equivalences being isomorphisms, and all morphisms being fibrations and
cofibrations, then our result applies and seems to be new even in this case: under
\Vopenkas principle any co-orthogonality class is generated by a set
of objects.

A counterexample to a similar question for orthogonality classes in locally presentable
categories was constracted in \cite{CS}. More specifically: it was
shown that there are orthogonality classes in the category of groups
which are not generated by any set of groups.

Additional motivation for the work on this paper was provided by Emmanuel
Farjoun who asked: whether or not any continuous, augmented, homotopy idempotent
functor from pointed simplicial sets to pointed simplicial sets is
equivalent, in a certain set-theoretical framework, to a
cellularization with respect to some simplicial set. This question is
motivated, in turn, by the following construction: given a nullification with
respect to a simplicial set $A$, the homotopy fibre $\bar P_A(X)$ of
the natural map $X \rarrow P_A(X)$ is a homotopy idempotent augmented
functor, which is usually not equivalent to $\CW_A X$. We prove in
this paper that there exist a simplicial set $B$ such that $\bar
P_A(X)$ is weakly equivalent to $\CW_B(X)$ for all $X$. An alternative
construction of $B$, not relying on \Vopenkas principle, was given in \cite{ChPS}.

\section{Simplicial orthogonality and co-orthogonality}
We will assume throughout the paper that \cat M is a simplicial model
category, but all the concepts in this section may be generalized to
arbitrary model categories using homotopy function complexes. For
every $A, X\in \cat M$ let $\tilde A$ denote the cofibrant replacement
of $A$ and let $\hat X$ denote the fibrant replacement of $X$.

We say that an object $X$ and a map $f\colon A\to B$ are
\textbf{homotopy orthogonal} or \textbf{simplicially orthogonal},
$f\hbot X$, if the induced map of simplicial sets 
\begin{equation}
f^*\colon \map(\tilde B,\hat X)\longrightarrow \map(\tilde A,\hat X)
\end{equation}
is a weak equivalence. 

Co-orthogonality is the dual concept: an object $A$ and a map $f\colon
X\rarrow Y$ are \textbf{homotopy co-orthogonal}, $A\htop f$ if the
induced map of simplicial sets 
\begin{equation}
f_*\colon \map(\tilde A,\hat X)\longrightarrow \map(\tilde A,\hat Y)
\end{equation}
is a weak equivalence. 

More generally, if $\calS$ is any class of maps, we denote by
${\calS}^{\hbot}$ the class of objects that are homotopy orthogonal to
all the maps in $\calS$, and say that objects in ${\calS}^{\hbot}$ are
\textbf{$\calS$-local}. The homotopy orthogonal complement of a class
$\calD$ of objects is defined similarly. Fibrant elements of
$f^{\hbot}$ are usually called $f$-local. 

For a class of morphisms $\calS$ in a model category $\cat M$, an \textbf{$\calS$-localization} is a homotopy endofunctor $L\colon \cat M \rarrow \cat M$ equipped with a natural transformation $\eta\colon \Id \rarrow L$ such that $\eta L \simeq L\eta$ and $L\eta\colon L\to LL$ is a weak equivalence on all objects, $\eta_X\colon X\to LX$ is in $({\calS}^{\hbot})^{\hbot}$ and $LX\in \calS^{\hbot}$ for all $X$ in \cat M. We also call it a \textbf{localization with respect to $\calS$} or, generically, a \textbf{homotopy localization}.

For a class of objects $\calA$ in a model category $\cat M$, an \textbf{$\calA$-co-localization} or \textbf{$\calA$-cellularization} is a homotopy endofunctor $C\colon \cat M \rarrow \cat M$ equipped with a natural transformation $\varepsilon\colon C \rarrow \Id$ such that $\varepsilon C \simeq C\varepsilon$ and $C\varepsilon \colon CC \to C$ is a weak equivalence on all objects, $\varepsilon_X\colon CX\to X$ is in ${\calA}^{\htop}$ and $CX\in (\calA^{\htop})^{\htop}$ for all $X$ in \cat M. We also call it a \textbf{co-localization with respect to $\calA$} or, generically, a \textbf{homotopy cellularization}.

Recall that a partially ordered set $P$ is called {\it $\lambda$-directed},
where $\lambda$ is a regular cardinal, if every subset of $P$ of cardinality
smaller than $\lambda$ has an upper bound.

\begin{lemma}\label{inj-fib}
Given a cofibrantly generated simplicial model category \cat M and a small category \cat C, consider the projective simplicial model structure on the category of functors $\cat M^{\cat C}$ described in \cite[Theorem~11.7.3]{Hirschhorn}. Suppose that $A$ is a cofibrant diagram in this model structure and $X$ is a fibrant object of \cat M, then the $\cat C^\op$-diagram $\hom(A,X)$ of simplicial sets is fibrant in the injective model structure on $\calS^{\cat C^\op}$.
\end{lemma}
\begin{proof}
We have to show that any commutative square
\[
\xymatrix{
C \ar@{^{(}->}[d]^{\dir{~}}_i \ar[r]  &  \hom(A,X) \ar@{->>}[d]\\
D \ar[r] & \ast,\\
}
\]
where $i$ is an injective (objectwise) trivial cofibration of $\cat C^\op$-diagrams of simplicial sets, admits a lift. By adjunction this problem is equivalent to finding a lift in the following commutative square in \cat M:
\[
\xymatrix{
A\otimes_{\cat C} C \ar[d] \ar[r]  &  X \ar@{->>}[d]\\
A\otimes_{\cat C} D \ar[r] & \ast.\\
}
\]
And this problem is equivalent, by another adjunction, to finding a lift in the following commutative square in $\cat M^\cat C$:
\[
\xymatrix{
\emptyset \ar@{^{(}->}[d] \ar[r]  &  X^D \ar@{->>}[d]^{\dir{~}}_{i^\ast}\\
A \ar[r] & X^C.\\
}
\]
In the last square the lift exists, since $A$ is projectively
cofibrant and $i^*$ is an objectwise trivial fibration, i.e., projective fibration.
\end{proof}

\begin{lemma}\label{cocomplete}
Let $\calS$ be any class of morphisms in a combinatorial simplicial model category $\cat M$,
and let $\calS^{\htop}=\calA$ be its homotopy co-orthogonal complement.
Then there exists a regular cardinal $\lambda$ such that $\calA$ is closed under 
$\lambda$-directed colimits in $\cat M$.
\end{lemma}

\begin{proof}
Let $I$ be a set of generating cofibrations for the model category $\cat M$. 
Choose a regular cardinal $\lambda$ such that any object of the set of domains and codomains of maps in $I$ 
is $\lambda$-presentable (such a cardinal exists since the category $\cat M$ is locally presentable). 
Let $P$ be any $\lambda$-directed partially ordered set, and suppose given
a diagram $A\colon P\to\cat M$. Let us depict it, for simplicity, as a chain:
\[
A_0 \longrightarrow A_1\longrightarrow \cdots\longrightarrow A_p\longrightarrow \cdots
\]
Suppose that the objects $A_p$ are in $\calA$ for each $p\in P$. We need to show that the colimit of this diagram is also in $\calA$.

Consider the category $\cat M^P$ of $P$-indexed diagrams in $\cat M$, and endow it with a model structure as described in \cite[11.6]{Hirschhorn}. Thus, weak equivalences and fibrations are objectwise, and cofibrations are retracts of free cell complexes. The diagram $A$ may be viewed as a single element in $\cat M^P$. Let $\tilde A$ be a cofibrant approximation of $A$ in the above model structure, hence obtaining the following commutative diagram in $\cat M$:
\[
\xymatrix{
\tilde A_0 \ar@{->>}[d]^{\dir{~}} \ar@{^{(}->}[r]  &  \tilde A_1 \ar@{->>}[d]^{\dir{~}} \ar@{^{(}->}[r] & \cdots \ar@{^{(}->}[r] & \tilde A_p \ar@{->>}[d]^{\dir{~}} \ar@{^{(}->}[r] & \cdots\\
A_0 \ar[r] & A_1 \ar[r] & \cdots \ar[r] & A_p \ar[r] & \cdots\\
}
\]
where $\tilde A_p$ is a cofibrant approximation of $A_p$ in \cat M. 

Note that $\colim \tilde A_p$ is a cofibrant object in \cat M, since $\colim \cat M^P \rarrow \cat M$ is a left Quillen functor \cite[11.6.8]{Hirschhorn}, hence preserves cofibrations.

For every map $f\colon X\rarrow Y$ in $\calS$, let $\hat f\colon \hat X\rarrow \hat Y$ be a fibrant approximation to $f$. The induced map
\[
\map(\colim \tilde A_p, \hat f)\colon \map(\colim \tilde A_p, \hat X) \longrightarrow \map(\colim \tilde A_p, \hat Y)
\]
can be written as
\[
\lim \map( \tilde A_p, \hat f)\colon \lim \map(\tilde A_p, \hat X) \longrightarrow \lim \map(\tilde A_p, \hat Y).
\]

By Lemma~\ref{inj-fib}, each of the $P^\op$-diagrams of simplicial sets $\map(\tilde A, \hat X)$ and $\map(\tilde A, \hat Y)$ is a fibrant object in the injective model structure on  the category of $P^\op$-diagrams of simplicial sets $\calS^{P^\op}$ , since $\tilde A$ is a cofibrant diagram in the projective model structure on $\cat M^{P^\op}$ and $\hat X, \hat Y$ are fibrant.
Therefore, $\map(\colim\tilde A_p,\hat f)$ is a homotopy inverse limit of weak equivalences,
so it is itself a weak equivalence. This shows that $\colim\tilde A_p$ is in $\calA$.

Trivial fibrations in $\cat M$ are preserved under $\lambda$-directed colimits, since the set of generating cofibrations has $\lambda$-presentable domains and codomains, hence there is a trivial fibration
\[
\xymatrix{
\colim \tilde A_p \ar@{->>}[r]^\sim  &  \colim A_p.
}
\]
We conclude that $\colim \tilde A_p$ is a cofibrant approximation of the $\colim A_p$, since any directed colimit of a projectively cofibrant diagram is cofibrant ($\tilde A$ is a cofibrant diagram in $\cat M^A$ and the colimit functor $\cat M^A\rarrow \cat M$ is left Quillen by \cite[11.6.8(1)]{Hirschhorn}). Hence, $\colim A_p$ is in $\calA$, as claimed.
\end{proof}

The next lemma significantly relies on \Vopenkas principle \cite{AR}.

\begin{lemma}\label{set-of-generators}
Suppose that \Vopenkas principle is true. Let $\calS$ be any class of morphisms in a 
combinatorial simplicial model category $\cat M$, and let ${\calD}={\calS}^{\htop}$. 
Then there exists a set of objects $\calX$ such that ${\calX}^{\htop} = {\calD}^{\htop}$.
\end{lemma}

\begin{proof}
By abuse of notation, we also denote by $\calD$ the full subcategory of $M$ generated by the class $\calD$. Since $\cat M$ is locally presentable, assuming \Vopenkas principle, it follows from \cite[Theorem~6.6]{AR} that $\calD$ is bounded, i.e., it has a small dense subcategory. We have shown in Lemma~\ref{cocomplete} that there exists a regular cardinal $\lambda$ such that $\calD$ is closed under $\lambda$-directed colimits. Hence, by \cite[Corollary~6.18]{AR}, the full subcategory generated by $\calD$ in $M$ is accessible. Thus, for a certain regular cardinal $\lambda_0 \geq \lambda$, the class $\calD$ contains a set $\calX$ of $\lambda_0$-presentable objects such that every object of $\calD$ is a $\lambda_0$-directed colimit of objects of $\calX$.

Since $\calX \subset \calD$, we have  $\calX^{\htop} \supset \calD^{\htop}$ and $(\calX^{\htop})^{\htop} 
\subset (\calD^{\htop})^{\htop} = \calS$. Our aim now is to show the reverse inclusion 
$(\calX^{\htop})^{\htop} \supset \calS$. By Lemma~\ref{cocomplete}, $(\calX^{\htop})^{\htop}$ is closed under $\lambda$-directed colimits .
Hence $(\calX^{\htop})^{\htop}$ is also closed under $\lambda_0$-directed colimits and  
every element of $\calD$ is a $\lambda_0$-directed colimit of elements of $\calX$. Then we can choose 
$\calX$ as our generating set.
\end{proof}

\begin{theorem}\label{HOSP}
Let $\cat M$ be a right proper, combinatorial, simplicial model category.
If \Vopenkas principle is assumed true, then
for any (possibly proper) class of objects $\calD$ there is a cellularization 
functor $CW_\calD$ with respect to $\calD$.
\end{theorem}
\begin{proof}
By Lemma~\ref{set-of-generators}, there exists a set $\calX$ of
objects in $\cat M$ such that
${\calX}^{{\htop}}={\calD}^{\htop}$. Then the cellularization with
respect to this set $\calX$ has ${\calD}^{\htop}$ as its class of
co-local equivalences, i.e. it is equivalent to the cellularization
with respect to $\calD$. Under the assumptions of the theorem, the
existence of the cellularization functor with respect to a set of
objects was established in \cite{Hirschhorn}.
\end{proof}

\section{Simplicial Idempotent functors and simplicial co-orthogonality}\label{Farjoun-problem}
In the proof of the main theorem bellow we rely on results, which are
dual to some of the statements in \cite{Farjoun:HHNC}. The assumption
of simpliciality (on functors) may be removed (if required) similarly to
\cite{cacho}.

\begin{theorem}
Let $\cat M$ be a cofibrantly generated simplicial model category. Let
$C$ be an augmented, homotopy idempotent simplicial functor from $\cat
M$ to $\cat M$. Then, assuming \Vopenkas principle, $C$ is equivalent
to homotopy cellularization with respect to some set of objects.
\end{theorem}

\begin{proof}
Formally, homotopy cellularization is a homotopy localization in the opposite category. The argument of \cite{Farjoun:HHNC} generalizes to any simplicial model category and continuous localization functor, since it does not use any small-object considerations. It shows that any simplicial localization functor is a localization with respect to the class $\calS$ of maps which are rendered into equivalences by the localization functor. Or, equivalently, this is a localization with respect to the class $\calD = \calS^{\hbot}$ of objects. In other words, the initial functor $C$ is a cellularization with respect to $\calD$.

By Lemma~\ref{set-of-generators} there exists a set $\calX$ of elements such that $\calD^{\htop}=\calX^{\htop}$. Therefore, the cellularization with respect to $\calX$ is equivalent to $C$.
\end{proof}

\bibliographystyle{abbrv}
\bibliography{Xbib}

\end{document}